\documentclass[12pt]{article} 
\usepackage{amssymb,amsmath,amsthm, eufrak, euscript}
\usepackage{epsfig}
\usepackage{graphicx}
\usepackage[margin=3cm]{geometry}

\begin{document}
\bibliographystyle{unsrt}

\newcommand{\im}{\Rightarrow}
\newcommand{\rr}{\ensuremath{\mathbb{R}}}
\newcommand{\cc}{\ensuremath{\mathbb{C}}}
\newcommand{\bs}{\backslash}
\newcommand{\summ}[3]{\ensuremath{\overset{#3}{\underset{#1=#2}{\sum}}}}
\newcommand{\summin}[2]{\ensuremath{\underset{#1 \in #2}{\sum}}}
\newcommand{\summunder}[1]{\ensuremath{\underset{#1}{\sum}}}
\newcommand{\cupp}[3]{\ensuremath{\overset{#3}{\underset{#1=#2}{\bigcup}}}}
\newcommand{\finvy}{\ensuremath{f^{-1}(y)}}
\newcommand{\ginvz}{\ensuremath{g^{-1}(z)}}
\newcommand{\dotcup}{\ensuremath{\overset{\cdot}{\cup}}}
\newenvironment{thm}[3][]{{\bf #2 Theorem #3} \em}{}
\newtheorem{lemma}{Lemma}
\newtheorem{claim}{Claim}
\newtheorem{prop}{Proposition}
\newtheorem{defn}{Definition}

\begin{center}
	{\large \bf Combinatorial Formulas for Classical Lie Weight Systems on Arrow Diagrams}
\end{center} 

\vspace{3mm}

\begin{center}
	Louis Leung\linebreak
	University of Toronto\linebreak
	louis.leung@utoronto.ca
\end{center}

\vspace{3mm}

\begin{abstract}
In \cite{haviv} Haviv gave a way of assigning Lie tensors to directed trivalent graphs.  Weight systems on oriented chord idagrams modulo 6T can then be constructed from such tensors.  In this paper we give explicit combinatorial formulas of weight systems using Manin triples constrcted from classical Lie algebras.  We then compose these oriented weight systems with the averaging map to get weight systems on unoriented chord diagrams and show that they are the same as the ones obtained by Bar-Natan (\cite{bn1}).  In the last section we carry out calculations on certain examples.   
\end{abstract}

\tableofcontents

\section{Introduction}

The study of finite type invariants of virtual knots (see \cite{gpv}) has led to interest in Gauss diagrams (circles with arrows joining distinct pairs of points on it) and arrow diagrams, which can be thought of as formal sums of Gauss diagrams modulo relations that correspond to the virtual Reidemeister moves.  It can be shown that a degree-$n$ finite type invariant gives rise to a functional on arrow diagrams with $n$ arrows.  in this paper we study such functionals coming from Lie bialgebras which are themselves constructed from the classical Lie algebras.  We begin with the following definitions.

\begin{defn}
Let $\Gamma$ be a disjoint union of oriented line segments and circles.  
\begin{enumerate}
	\item A chord diagram with skeleton $\Gamma$ is a diagram with chords joining distinct pairs of points on $\Gamma$. 
	\item An oriented chord diagram is a chord diagram together with an orientation on each of its chords. 
	\item A trivalent diagram with skeleton $\Gamma$ is $\Gamma$ together with a graph $G$ with univalent and trivalent vertices so that all the univalent vertices are attached to distinct points on $\Gamma$, and each trivalent vertex comes with an orientation (i.e. an ordering of the incident edges modulo cyclic permutations).
	\item An arrow diagram is a trivalent diagram with an orientation on each of the edges of $G$.
	\item An arrow diagram is acyclic if the underlying graph $G$ contains no oriented cycles.
\end{enumerate}
\end{defn}

\begin{figure}[h!]
	\centerline{\epsfig{file=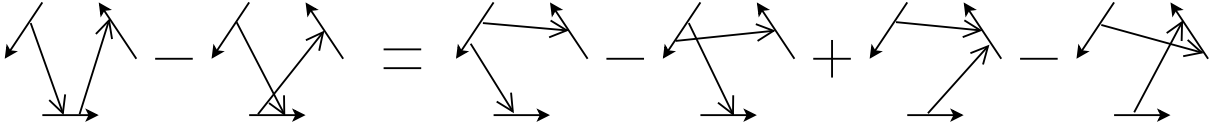, scale=0.25}}
	\caption{The 6T relation.}
      \label{fig:6t}
\end{figure}

We consider the vector space generated by oriented chord diagrams modulo the relation 6T.  (For origins see \cite{gpv}.)  It is depicted in figure \ref{fig:6t}.

\begin{figure}[h!]
	\centerline{\epsfig{file=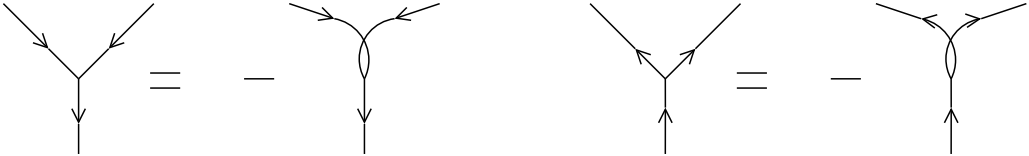, scale=0.20}}
	\caption{AS: the orientations of vertices change at the cost of a minus sign.}
      \label{fig:as}
\end{figure}

\begin{figure}[h!]
	\centerline{\epsfig{file=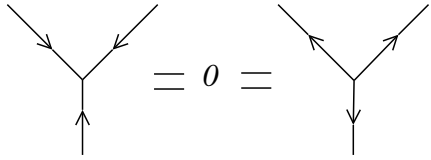, scale=0.20}}
	\caption{NS: There are no sinks (left) and no sources (right).}
      \label{fig:ns}
\end{figure}

\begin{figure}[h!]
	\centerline{\epsfig{file=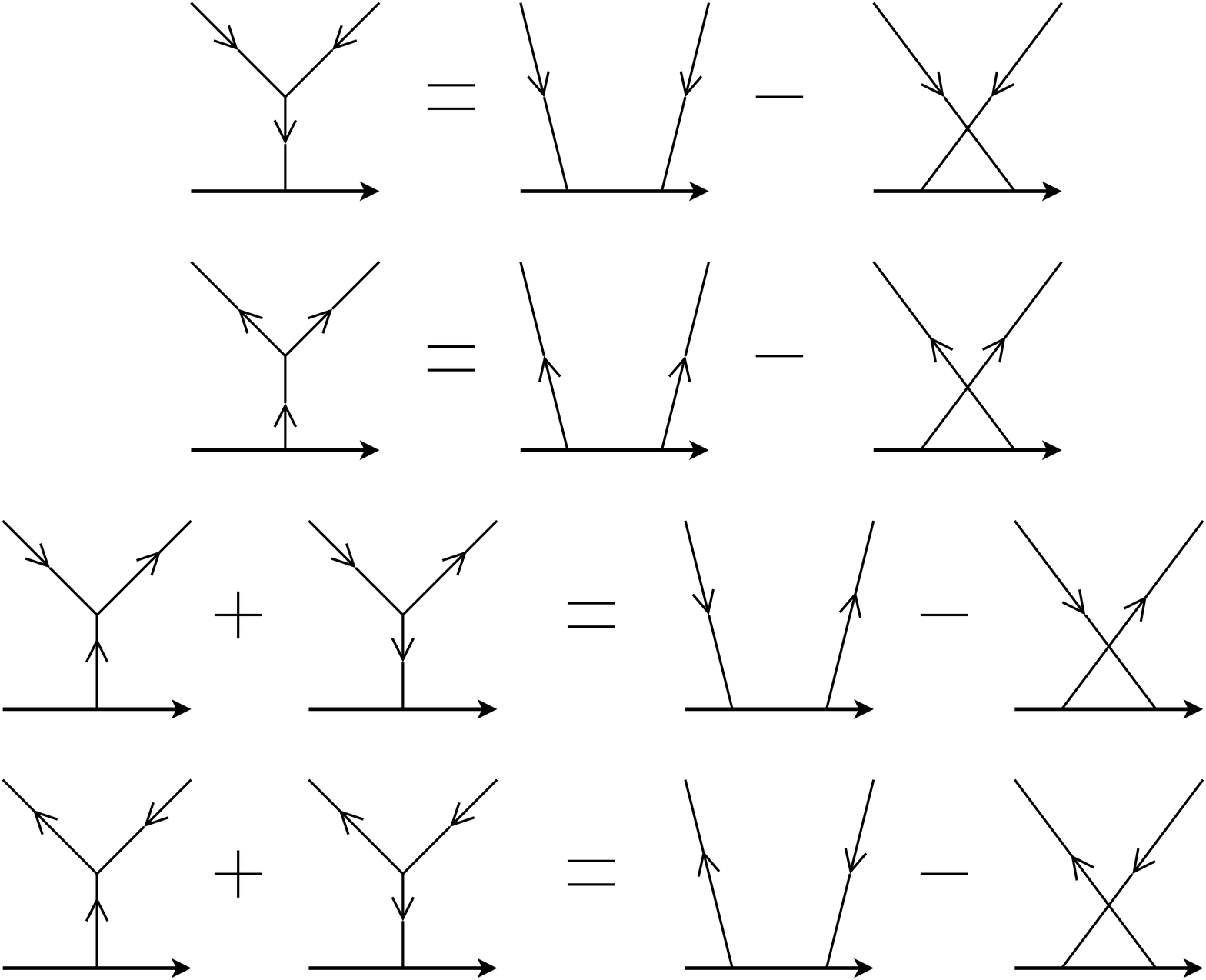, scale=0.20}}
	\caption{The STU relations}
      \label{fig:stu}
\end{figure}

\begin{figure}[h!]
	\centerline{\epsfig{file=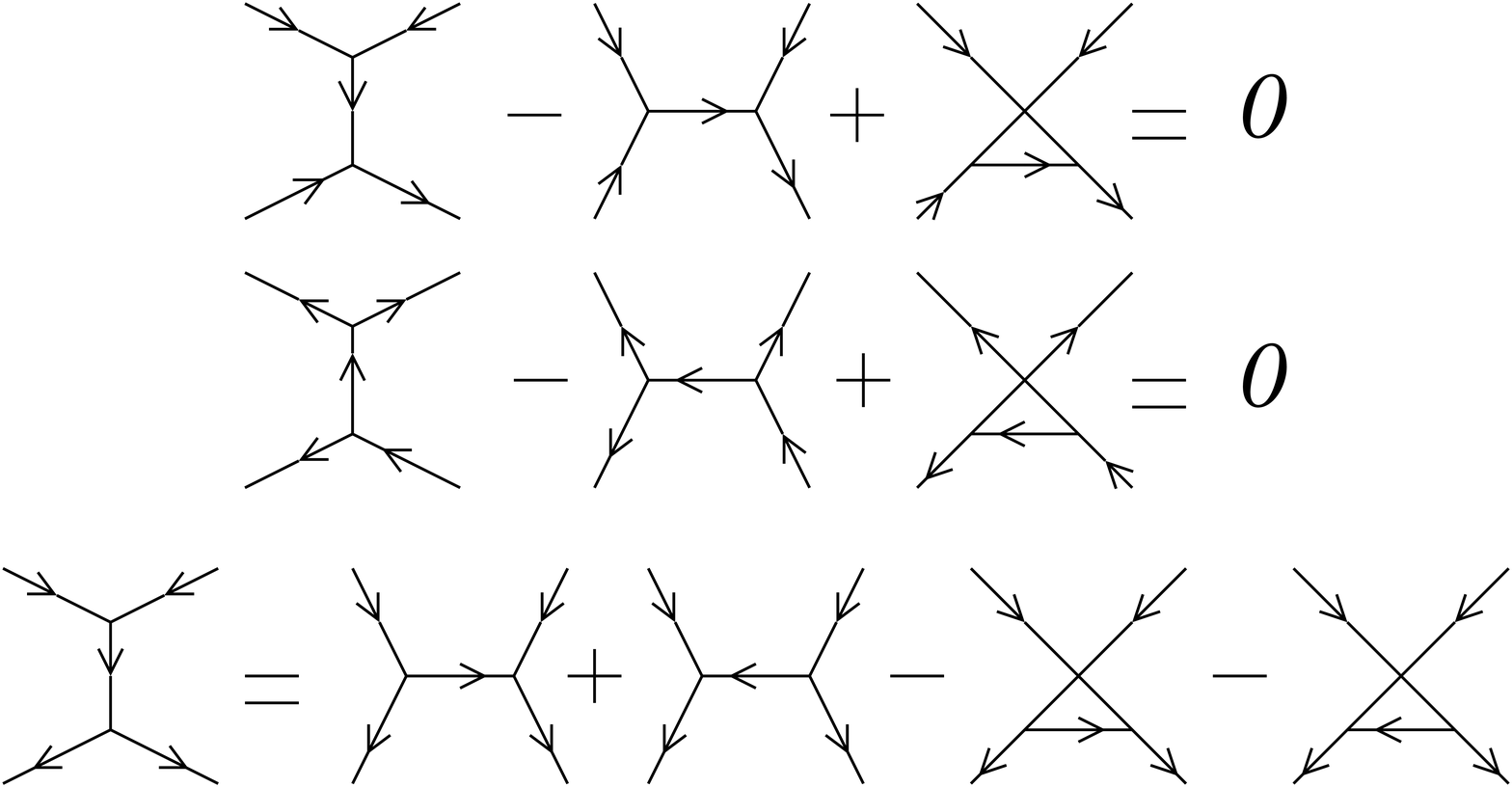, scale=0.20}}
	\caption{The IHX relations.}
      \label{fig:ihx}
\end{figure}

Similarly we consider the vector space generated by arrow diagrams modulo the relations $AS$ (antisymmetry), $NS$ (no sink, no source), $STU$ and $IHX$ in the space of arrow diagrams.  (See figures \ref{fig:as} to \ref{fig:ihx}.)

\begin{defn}
$\vec{\mathcal{C}}(\Gamma)$ is the vector space of oriented chord diagrams on $\Gamma$ modulo 6T.  $\vec{\mathcal{A}}(\Gamma)$ is the vector space of arrow diagrams on $\Gamma$ modulo $AS$, $NS$, $STU$ and $IHX$.  From now on we will use `oriented chord diagrams' and `arrow diagrams' to refer to equivalent classes in $\vec{\mathcal{C}}(\Gamma)$ and $\vec{\mathcal{A}}(\Gamma)$, respectively.  A functional from $\vec{\mathcal{C}}(\Gamma)$ or $\vec{\mathcal{A}}(\Gamma)$ to $\mathbb{C}$ is called a weight system.
\end{defn}

Our goal is to study weight systems on $\vec{\mathcal{C}}(\Gamma)$.  Polyak (\cite{polyak}) proved that the space of oriented chord diagrams modulo 6T is isomorphic to the subspace of $\vec{\mathcal{A}}(\Gamma)$ generated by acyclic arrow diagrams.  We will first study functionals on $\vec{\mathcal{A}}(\Gamma)$, which are easier to construct since they have a close realtion to Lie bialgebras.  Given Lie bialgebras constructed from classical Lie algebras we will then give algorithms to compute the corresponding weight systems on diagrams in $\vec{\mathcal{C}}(\Gamma)$

In section \ref{sec:bialg} we will follow \cite{es} to construct a Manin triple from a simple Lie algebra.  In section \ref{sec:tensor} we follow \cite{haviv} to assign tensors to elements of $\vec{\mathcal{A}}(\Gamma)$ and construct weight systems when the skeleton is a circle.  In section \ref{sec:calc}, using Manin triples coming from the families $\mathfrak{gl}$, $\mathfrak{so}$ and $\mathfrak{sp}$, we give explicit formulas to turn $\vec{\mathcal{C}}(\Gamma)$ into tensors.  We also introduce the averaging map and compose it with our tensors and compare the results with those in \cite{bn1} on unoriented diagrams.  Finally in section \ref{sec:sample} we apply the results from section \ref{sec:calc} and do some sample calculations on $\vec{\mathcal{C}}(S^{1})$.

\subsection{Acknowledgement}  This paper is part of the author's Ph.D. research at the University of Toronto under the supervision of Dror Bar-Natan.  The author would like to thank him for his guidance and suggestions in the writing of this paper.

\section{Lie bialgebras and Manin triples from a simple Lie algebra}
\label{sec:bialg}

\begin{defn}
A Lie bialgebra $(\mathfrak{g},[,],\delta )$ is a Lie algebra $(\mathfrak{g},[,])$ equipped with an antisymmetric $cobracket$ map $\delta :\mathfrak{g}\rightarrow\mathfrak{g}\otimes\mathfrak{g}$ satisfying the coJacobi identity
\begin{equation}
\nonumber (id+\tau +\tau^{2})((\delta\otimes id)\delta (x))=0
\end{equation}
and the cocycle condition
\begin{equation}
\nonumber \delta([xy])=ad_{x}(\delta y)-ad_{y}(\delta x),
\end{equation}
for any $x,y\in\mathfrak{g}$, where $\tau$ is the cyclic permutation on $\mathfrak{g}^{\otimes 3}$. 
\end{defn}

\begin{defn}
A finite dimensional Manin triple is a triple of finite dimensional Lie algebras $(\tilde{\mathfrak{g}},\mathfrak{g}_{+},\mathfrak{g}_{-})$, where $\tilde{\mathfrak{g}}$ is equipped with a metric (a symmetric nondegenerate invariant bilinear form) $(.,.)$ such that
\begin{enumerate}
	\item $\tilde{\mathfrak{g}}=\mathfrak{g}_{+}\oplus\mathfrak{g}_{-}$ as a vector space and $\mathfrak{g}_{+},\mathfrak{g}_{-}$ are Lie subalgebras of $\tilde{\mathfrak{g}}$.
	\item $\mathfrak{g}_{+},\mathfrak{g}_{-}$ are isotropic with respect to $(.,.)$.
\end{enumerate}
\end{defn}

As a consequence, $\mathfrak{g}_{+}$ and $\mathfrak{g}_{-}$ are maximal isotropic subalgebras.  Suppose $(\tilde{\mathfrak{g}},\mathfrak{g}_{+},\mathfrak{g}_{-})$ is a Manin triple.  The metric then induces a nondegenerate pairing $\mathfrak{g}_{+}\otimes\mathfrak{g}_{-}\rightarrow\mathbb{C}$, and hence a Lie algebra structure on $\mathfrak{g}_{+}^{*}\cong\mathfrak{g}_{-}$.  Let $\delta$ be the induced coalgebra structure on $\mathfrak{g}_{+}$.  We can check by direct computation (\cite{es}) that the cocycle condition is satisfied.  $(\mathfrak{g}_{+},[.,.],\delta)$ is therefore a Lie bialgebra. 

In fact the process can be reversed.  Given a Lie bialgebra $\mathfrak{g}$, we may define a symmetric nondegenerate bilinear form $(.,.)_{\mathfrak{g}\oplus\mathfrak{g}^{*}}$ by $((x,f),(x',f'))_{\mathfrak{g}\oplus\mathfrak{g}^{*}}=f(x')+f'(x)$.  If ${e_{i}}$ is a basis of $\mathfrak{g}$ and ${f^{i}}$ is a basis of $\mathfrak{g}^{*}$ with $[e_{i},e_{j}]=c_{ij}^{k}e_{k}$ and $[f^{r},f^{s}]=\gamma^{rs}_{t}f^{t}$, then we can define a Lie algebra structure on $\mathfrak{g}\oplus\mathfrak{g}^{*}$ by

\begin{equation}
[f^{r},e_{s}]=c^{r}_{st}f^{t}-\gamma^{rt}_{s}e_{t}
\label{eqn:bracket}
\end{equation}

\noindent and keeping the bracket for $\mathfrak{g}$ and $\mathfrak{g}^{*}$, with the consequence that $(.,.)_{\mathfrak{g}\oplus\mathfrak{g}^{*}}$ is invariant.  (See section 1.3 of \cite{cp}.) $(.,.)_{\mathfrak{g}\oplus\mathfrak{g}^{*}}$ is therefore a metric and $(\mathfrak{g}\oplus\mathfrak{g}^{*},\mathfrak{g},\mathfrak{g}^{*})$ is a Manin triple. 

There is a standard way to obtain Manin triples from simple Lie algebras, and those are the ones we are going to use.  The construction below follows Chapter 4 of \cite{es}.  Given a simple Lie algebra over $\mathbb{C}$ with metric $(.,.)$, we fix a Cartan subalgebra $\mathfrak{h}$ and consider the corresponding positive and negative root spaces $\mathfrak{n}_{+}$ and $\mathfrak{n}_{-}$.  For each root $\alpha$ we consider $e_{\alpha}\in\mathfrak{g}_{\alpha}$ and $f_{\alpha}\in\mathfrak{g}_{-\alpha}$ (where $\mathfrak{g}_{\pm\alpha}$ are the root spaces corresponding to $\pm\alpha$) such that $(e_{\alpha},f_{\alpha})=1$.  Let $h_{\alpha}=[e_{\alpha},f_{\alpha}]$.  We consider the Lie algebra

\begin{equation}
\nonumber \tilde{\mathfrak{g}}=\mathfrak{n}_{+}\oplus\mathfrak{h}^{(1)}\oplus\mathfrak{h}^{(2)}\oplus\mathfrak{n}_{-}
\end{equation}

\noindent where $\mathfrak{h}^{(1)}\cong\mathfrak{h}\cong\mathfrak{h}^{(2)}$ and with bracket defined by:

\begin{center}
$[\mathfrak{h}^{(1)},\mathfrak{h}^{(2)}]=0$,\hspace{7mm} $[\mathfrak{h}^{(i)},e_{\alpha}]=\alpha (h)e_{\alpha}$,\\ $[\mathfrak{h}^{(i)},f_{\alpha}]=-\alpha (h)f_{\alpha}$,\hspace{5mm} and \hspace{5mm} $[e_{\alpha},f_{\alpha}]=\frac{1}{2}(h^{(1)}_{\alpha}+h^{(2)}_{\alpha})$.
\end{center}

We define the following metric on $\tilde{\mathfrak{g}}$:

\begin{equation}
\nonumber (x+h^{(1)}+h^{(2)},x'+h'^{(1)}+h'^{(2)})_{\tilde{\mathfrak{g}}}=2((h^{(1)},h'^{(2)})_{\mathfrak{g}}+(h'^{(1)},h^{(2)})_{\mathfrak{g}})+(x,x')_{\mathfrak{g}}
\end{equation}

We can check that $(\tilde{\mathfrak{g}},\mathfrak{n}_{+}\oplus\mathfrak{h}^{(1)},\mathfrak{n}_{-}\oplus\mathfrak{h}^{(2)})$ is a Manin triple.  In fact $\tilde{\mathfrak{g}}$ is a Lie bialgebra with r-matrix

\begin{equation}
\nonumber \tilde{r}=\sum_{\alpha\in\Delta^{+}}e_{\alpha}\otimes f_{\alpha}+\frac{1}{2}\sum_{i}k_{i}^{(1)}\otimes k_{i}^{(2)},
\end{equation}

\noindent i.e., $\delta(x)=ad_{x}(\tilde{r})$, where $\{k_{i}\}$ is an orthonormal basis of $\mathfrak{h}$ with respect to $(.,.)$.  We define the projection $\pi :\tilde{\mathfrak{g}}\rightarrow\mathfrak{g}$ where

\begin{center}
$\pi|_{\mathfrak{n}_{+}\oplus\mathfrak{n}_{-}}=Id$\hspace{7mm}$\pi(h^{(1)}_{\alpha})=h_{\alpha}=\pi(h^{(2)}_{\alpha})$
\end{center}

\noindent This map endows $\mathfrak{g}$ with a quasitriangular Lie bialgebra structure with r matrix

\begin{equation}
\nonumber r=\sum_{\alpha\in\Delta^{+}}e_{\alpha}\otimes f_{\alpha}+\frac{1}{2}\sum_{i}k_{i}\otimes k_{i},
\end{equation}

\noindent so $\delta(x)=ad_{x}(r)$.  The Lie subalgebras $\mathfrak{b}_{+}=\mathfrak{n}_{+}\oplus\mathfrak{h}$ and $\mathfrak{b}_{-}=\mathfrak{n}_{-}\oplus\mathfrak{h}$ are Lie subbialgebras.

Note that the map $\pi$ is a Lie algebra homomorphism.  If $\mathfrak{g}$ is given as a matrix Lie algebra then $\pi$ is a representation.

\section{Directed trivalent graphs and Lie tensors}
\label{sec:tensor}

Since Polyak (\cite{polyak}) has proved that $\vec{\mathcal{C}(\Gamma)}$ modulo 6T is isomorphic to subspace of $\vec{\mathcal{A}}(\Gamma)$ generated by the acyclic diagrams, a functional on $\vec{\mathcal{A}}(\Gamma)$ induces a functional on $\vec{\mathcal{C}(\Gamma)}$.  In this section we construct a funtional on $\vec{\mathcal{A}}(\Gamma)$ by following Haviv's method of assigning Lie tensors to directed trivalent graphs and a representation to the skeleton (\cite{haviv}) and then calculate the trace.  Let $\mathfrak{g}$ be a Lie bialagebra with basis \{$e_{i}$\} and let \{$f^{i}$\} be the dual basis of $\mathfrak{g}^{*}$.  We would like to assign to each directed trivalent graph a tensor in $\tilde{\mathfrak{g}}=\mathfrak{g}\oplus\mathfrak{g}^{*}$.  First to the straight arrow we assign the element $f^{i}\otimes e_{i}$, which coresponds to the identity map.

\begin{figure}[h]
	\centerline{\epsfig{file=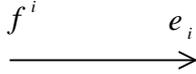, scale=0.25}}
	\caption{The arrow.}
      \label{fig:arrowtensor}
\end{figure}

\begin{figure}[h]
	\centerline{\epsfig{file=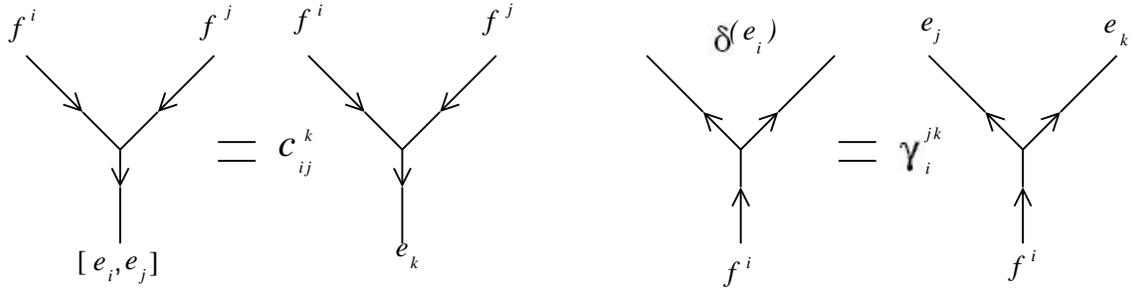, scale=0.25}}
	\caption{The bracket (two in, one out) tensor (left) and the cobracket(one in, two out) tensor (right).}
      \label{fig:bracket}
\end{figure}

Given the NS relation we have two types of vertices (one in, two out and two in, one out).  For the first type we assign to it the cobracket tensor $f^{i}\otimes\delta(e_{i})=\gamma^{jk}_{i}(f^{i}\otimes e_{j}\otimes e_{k})$, while to the second we assign the bracket tensor $f^{i}\otimes f^{j}\otimes [e_{i},e_{j}]=c^{k}_{ij}(f^{i}\otimes f^{j}\otimes e_{k})$, where $c^{k}_{ij}$ and $\gamma^{jk}_{i}$ are the structure constants for the bracket and the cobracket, respectively.  (See figure \ref{fig:bracket}.)  It is worth noting that under this interpretation, the STU becomes the diagrammatic way of saying $[x,y]=xy-yx$ in the universal enveloping algebra of $\tilde{\mathfrak{g}}$ (see equation \ref{eqn:bracket}), the 3-term IHX relations become the Jacobi and coJacobi identities, and the 5-term IHX becomes the cocycle identity.

To assign a tensor to a directed trivalent graph, we break it down into subgraphs with 0 or 1 vertex and at the points of gluing we contract using the metric.

\begin{figure}[h!]
	\centerline{\epsfig{file=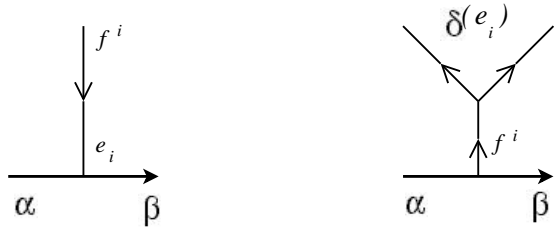, scale=0.25}}
	\caption{We assign a representation to the skeleton.  The diagrams correspond to $(R(e_{i}))^{\beta}_{\alpha}f^{i}$ and $(R(f^{i}))^{\beta}_{\alpha}(\delta{e_{i}})$, respectively.}
      \label{fig:repexample}
\end{figure}

Given a representation $R:\tilde{\mathfrak{g}}\rightarrow End(V)$ where $V$ has a basis $b=\{v_{1},...,v_{d}\}$, if the skeleton is part of the picture, we assign Greek letters ranging over $\{1,...,d\}$ to each section of the skeleton.  For example in figure \ref{fig:repexample}, we assign $v^{\beta}(R(e_{i})(v_{\alpha}))f^{i}=(R(e_{i}))^{\beta}_{\alpha}f^{i}$ and $v^{\beta}(R(f^{i})(v_{\alpha}))(\delta{e_{i}})=(R(f^{i}))^{\beta}_{\alpha}(\delta{e_{i}})\in\tilde{\mathfrak{g}}\otimes\tilde{\mathfrak{g}}$, respectively, to the pictures.  Note, since $v^{\beta}=\langle v_{\beta},.\rangle$ where $\langle .,.\rangle$ is the inner product of $\mathbb{C}^{d}$ with respect to the given basis, the same values can be written as $\langle v_{\beta},R(e_{i})(v_{\alpha})\rangle f^{i}$ and $\langle v_{\beta},R(f^{i})(v_{\alpha})\rangle\delta{e_{i}}$.  Here we do not distinguish $R(f^{i})$ or $R(e_{i})$ from its matrix with respect to the basis $b$.

If we restrict ourselves to the case where the skeleton is a circle, we can see that the construction above gives us the trace of the tensor in the given representation.

\section{Calculations using Lie bialgebras coming from classical Lie algebras}
\label{sec:calc}

In this section we will calculate weight systems coming from Manin triples constructed from classical Lie algebras.  The metric in each of the Lie algebras is $(A,B)=tr(AB)$. Throughout this section $e_{ij}$ or $e^{ij}$ is the matrix whose $ij$-th entry is 1 and zero everywhere else.  Let $\pi$ be the map given in section \ref{sec:bialg}.  We think of it as a representation and let $x_{i}=\pi(e_{i})$ and $\xi^{i}=\pi(f^{i})$.

\subsection{$\mathfrak{gl}(N)$}

We begin with $\mathfrak{gl}(N)$.  Note that $\mathfrak{gl}(N)$ is not simple, but our construction in section \ref{sec:bialg}, when applied to $\mathfrak{sl}(N)$, results in the Manin triple $(\tilde{\mathfrak{u}},\mathfrak{u},\mathfrak{l})$, where $\mathfrak{u}$ is the Lie algebra of (non-strictly) upper triangular matrices with trace 0 and its dual $\mathfrak{l}$ is the lower triagular matrices with trace 0.  This implies that both $\mathfrak{u}$ and $\mathfrak{l}$ are Lie bialgebras.  Let $\mathfrak{s}$ be the commutative Lie algebra of scalar matrices,  then $\mathfrak{u}\oplus\mathfrak{s}$ is the Lie algebra of upper triangular matrices.  In fact it is a Lie bialgebra if we set $\delta(y)=0$ for any $y\in\mathfrak{s}$.  Now following section \ref{sec:bialg} we consider $\mathfrak{g'}=\mathfrak{u}\oplus\mathfrak{s}^{(1)}\oplus\mathfrak{s}^{(2)}\oplus\mathfrak{l}$ where $\mathfrak{s}^{(1)}\cong\mathfrak{s}\cong\mathfrak{s}^{(2)}$ and consider the metric $(.,.)'$ and the bracket $[.,.]'$ on $\mathfrak{g'}$ defined as follows:

\begin{equation}
\nonumber (x+d^{(1)}+d^{(2)},x'+d'^{(1)}+d'^{(2)})'=2((d^{(1)},d'^{(2)})+(d'^{(1)},d^{(2)}))+(x,x')_{\tilde{\mathfrak{u}}}
\end{equation}

\begin{center}
$[\mathfrak{u},\mathfrak{l}]'=[\mathfrak{u},\mathfrak{l}]_{\tilde{\mathfrak{u}}}$\hspace{7mm}$[\mathfrak{s}^{(i)},x]=0$ for $i=1,2$ and any $x\in\mathfrak{g'}$,
\end{center}

\noindent where $(.,.)_{\tilde{\mathfrak{u}}}$ and $[.,.]_{\tilde{\mathfrak{u}}}$ are the metric and the bracket in $\tilde{\mathfrak{u}}$, respectively.

The dual to $\mathfrak{u}\oplus\mathfrak{s}^{(1)}$ under $(.,.)'$ is then $\mathfrak{l}\oplus\mathfrak{s}^{(2)}$.  If we let $\mathfrak{g}=\mathfrak{u}\oplus\mathfrak{s}$ and $\mathfrak{g'}=\mathfrak{g}\oplus\mathfrak{g}^{*}$ then we consider $(\mathfrak{g'},\mathfrak{g},\mathfrak{g}^{*})$ to be the Manin triple constructed from $\mathfrak{gl}(N)$.  In this paper we use $\mathfrak{gl}(N)$ instead of $\mathfrak{sl}(N)$.  The map $\pi$ is defined similarly as in section \ref{sec:bialg}.

\begin{figure}[h]
	\centerline{\epsfig{file=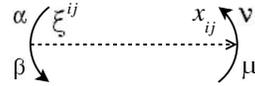, scale=0.25}}
	\caption{The simplest subdiagram in an arrow diagram.}
      \label{fig:glnonearrow}
\end{figure}

We consider the tensor in figure \ref{fig:glnonearrow}.

Note that $\{e_{ij}\}_{i\leq j}$ forms a basis of $\mathfrak{g}$ and $\{f^{ij}\}_{i\leq j}$ (where $f^{ij}=(\frac{1}{2})^{\delta^{ij}}e^{ji}$) forms the corresponding dual basis of $\mathfrak{g}^{*}$ under $(.,.)$. Using Haviv's method we assign $f^{ij}\otimes e_{ij}$ to each chord (assuming the Einstein summation convention).  We have $x_{ij}=e_{ij}$ and $\xi^{ij}=(\frac{1}{2})^{\delta^{ij}}e^{ji}$ (so $e$ and $x$, and similarly $f$ and $\xi$, are the same matrices sitting in different spaces).  The tensor is therefore given by 

\begin{eqnarray}
\nonumber 
T_{\alpha\mu}^{\beta\nu} =\sum_{i<j}(\xi^{ij})^{\beta}_{\alpha}(x_{ij})^{\nu}_{\mu}+\frac{1}{2}\sum_{i=j}(\xi^{ij})^{\beta}_{\alpha}(x_{ij})^{\nu}_{\mu}\\
\nonumber
=\sum_{i<j}\delta^{j\beta}\delta^{i}_{\alpha}\delta_{i}^{\nu}\delta_{j\mu}+\frac{1}{2}\sum_{i=j}\delta^{j\beta}\delta^{i}_{\alpha}\delta_{i}^{\nu}\delta_{j\mu}
\end{eqnarray}

\begin{center}
$=\left\{\begin{array}{ll}
\delta^{\beta}_{\mu}\delta_{\alpha}^{\nu} & \textnormal{for $\alpha < \beta$}\\
\frac{1}{2}\delta^{\beta}_{\mu}\delta_{\alpha}^{\nu} & \textnormal{for $\alpha = \beta$}\\
0 & \textnormal{otherwise.}
\end{array}\right.$
\end{center}

\begin{figure}[h]
	\centerline{\epsfig{file=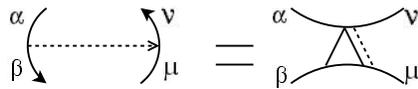, scale=0.25}}
	\caption{A diagrammatic representation of the $\mathfrak{gl}(N)$ tensor.}
      \label{fig:gln}
\end{figure}

The result can be expressed diagramatically as in figure \ref{fig:gln}, where the dotted equality sign means that when $\alpha = \beta =\mu =\nu$, the resulting value is $\frac{1}{2}$.

\subsection{$\mathfrak{so}(2N)$}

For $\mathfrak{so}(2N)$ we use the basis $x_{ijk}$ and $\xi^{ijk}$ where $0\leq i\leq j\leq N$ and $k=1,2$.  We have: 

\begin{eqnarray}
\nonumber x_{ij1}&=&e_{ij}-e_{j+N,i+N}\\
\nonumber x_{ij2}&=&e_{i,j+N}-e_{j,i+N}\\
\nonumber \xi^{ij1}&=&(\frac{1}{2})^{\delta^{ij}+1}(e^{ji}-e^{i+N,j+N})\\
\nonumber \xi^{ij2}&=&(\frac{1}{2})^{\delta^{ij}+1}(e^{j+N,i}-e^{i+N,j})
\end{eqnarray}

Note this is not exactly a basis because $x_{ii2}=\xi^{ii2}=0$ but that only means in the calculation of $T_{\alpha\mu}^{\beta\nu}$ we have some extra zero terms.  $T_{\alpha\mu}^{\beta\nu}$ in this case is given by

\begin{eqnarray}
\nonumber T_{\alpha\mu}^{\beta\nu}&=&\sum_{i\leq j}(\xi^{ijk})_{\alpha}^{\beta}(x_{ijk})_{\mu}^{\nu}\\
\nonumber &=&\sum_{i\leq j}((\frac{1}{2})^{\delta^{ij}+1}(\delta^{\beta,j}\delta^{i}_{\alpha}-\delta^{\beta,i+N}\delta^{j+N}_{\alpha})(\delta^{\nu}_{i}\delta_{j\mu}-\delta^{\nu}_{j+N}\delta_{i+N,\mu})\\ 
\nonumber &&+(\frac{1}{2})^{\delta^{ij}+1}(-\delta^{\beta,i+N}\delta^{j}_{\alpha}+\delta^{\beta,j+N}\delta^{i}_{\alpha})(\delta^{\nu}_{i}\delta_{j+N,\mu}-\delta^{\nu}_{j}\delta_{i+N,\mu}))
\end{eqnarray}

\begin{figure}[h]
	\centerline{\epsfig{file=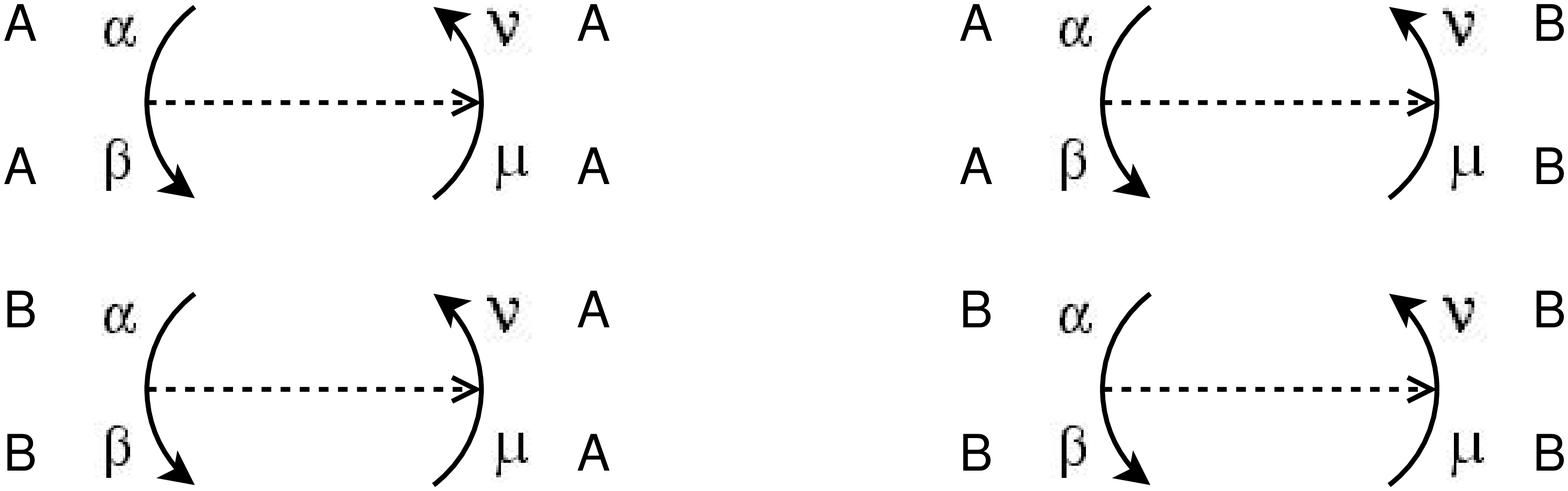, scale=0.25}}
	\caption{Four ways to assign $A$ and $B$ to $\alpha$, $\beta$, $\mu$ and $\nu$.}
      \label{fig:so2nABk1}
\end{figure}

In the expression above the first product is given by the case when $k=1$ while the second product corresponds to the case $k=2$.  We consider the first product.  Letting $A$ be the set $\{1,2,...,N\}$ and $B$ the set $\{N+1,N+2,...,2N\}$, we can assign $A$ or $B$ to each of $\alpha$, $\beta$, $\mu$ and $\nu$ to indicate which set it belongs to.  The four combinations that correspond to the four ways to expand the first summand in the last equation are as given in figure \ref{fig:so2nABk1}.

\begin{figure}[h!]
	\centerline{\epsfig{file=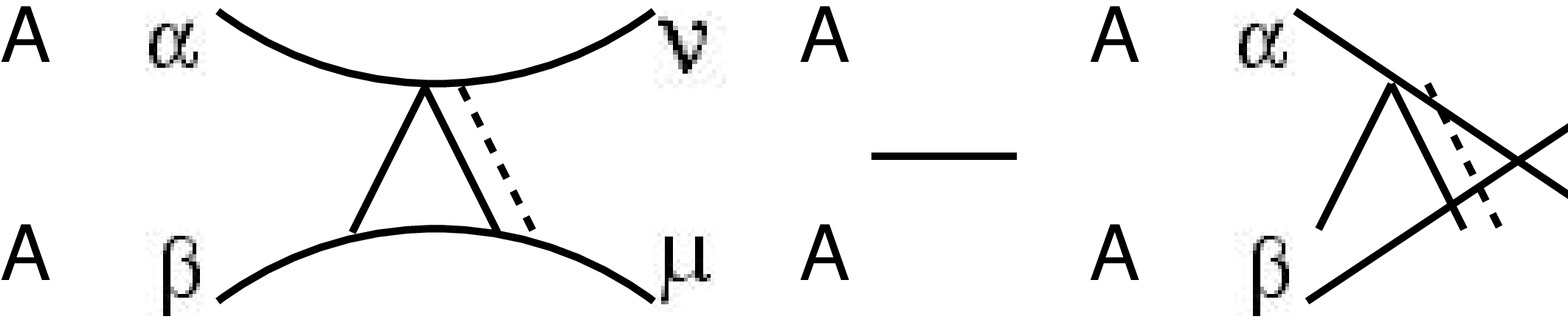, scale=0.25}}
	\caption{The $\mathfrak{so}(2N)$ tensor when $k=1$.}
      \label{fig:so2nk1}
\end{figure}

\noindent Starting at the top left corner and going clockwise the four pictures correspond to $\delta^{\beta,j}\delta^{i}_{\alpha}\delta^{\nu}_{i}\delta_{j\mu}$, $\delta^{\beta,j}\delta^{i}_{\alpha}\delta^{\nu}_{j+N}\delta_{i+N,\mu}$, $\delta^{\beta,i+N}\delta^{j+N}_{\alpha}\delta^{\nu}_{j+N}\delta_{i+N,\mu}$ and $\delta^{\beta,i+N}\delta^{j+N}_{\alpha}\delta^{\nu}_{i}\delta_{j\mu}$, respectively. If we delete the arrow and join two Greek letters to indicate that they are equal (or equal modulo $N$, if one comes from $A$ and the other from $B$), then the first summand can be expressed diagramatically as in figure \ref{fig:so2nk1}.

\begin{figure}[h!]
	\centerline{\epsfig{file=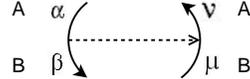, scale=0.25}}
	\caption{The only possible assignment of $A$ and $B$ when $k=2$.}
      \label{fig:so2nABk2}
\end{figure}

\begin{figure}[h!]
	\centerline{\epsfig{file=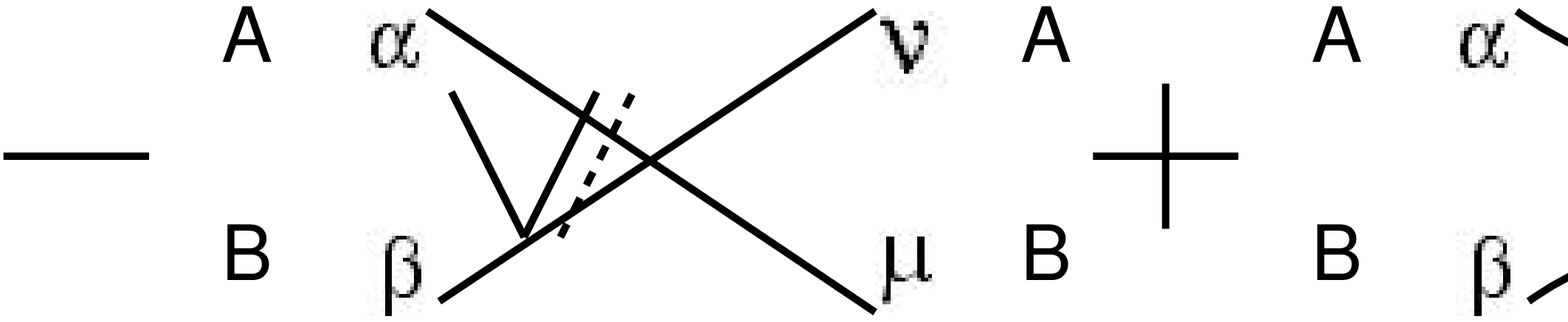, scale=0.25}}
	\caption{The $\mathfrak{so}(2N)$ tensor when $k=2$.}
      \label{fig:so2nk2}
\end{figure}

The inequality comes from the fact that $i\leq j$ and again the dotted equality means that when $\alpha=\beta=\mu=\nu$ the resulting value comes with an extra factor of $\frac{1}{2}$.  Now we look at the product corresponding to the case $k=2$.  The only possible assignment is given in figure \ref{fig:so2nABk2}.  Similarly we can represent the second summand diagramatically.  See figure \ref{fig:so2nk2}.

\begin{figure}[h!]
	\centerline{\epsfig{file=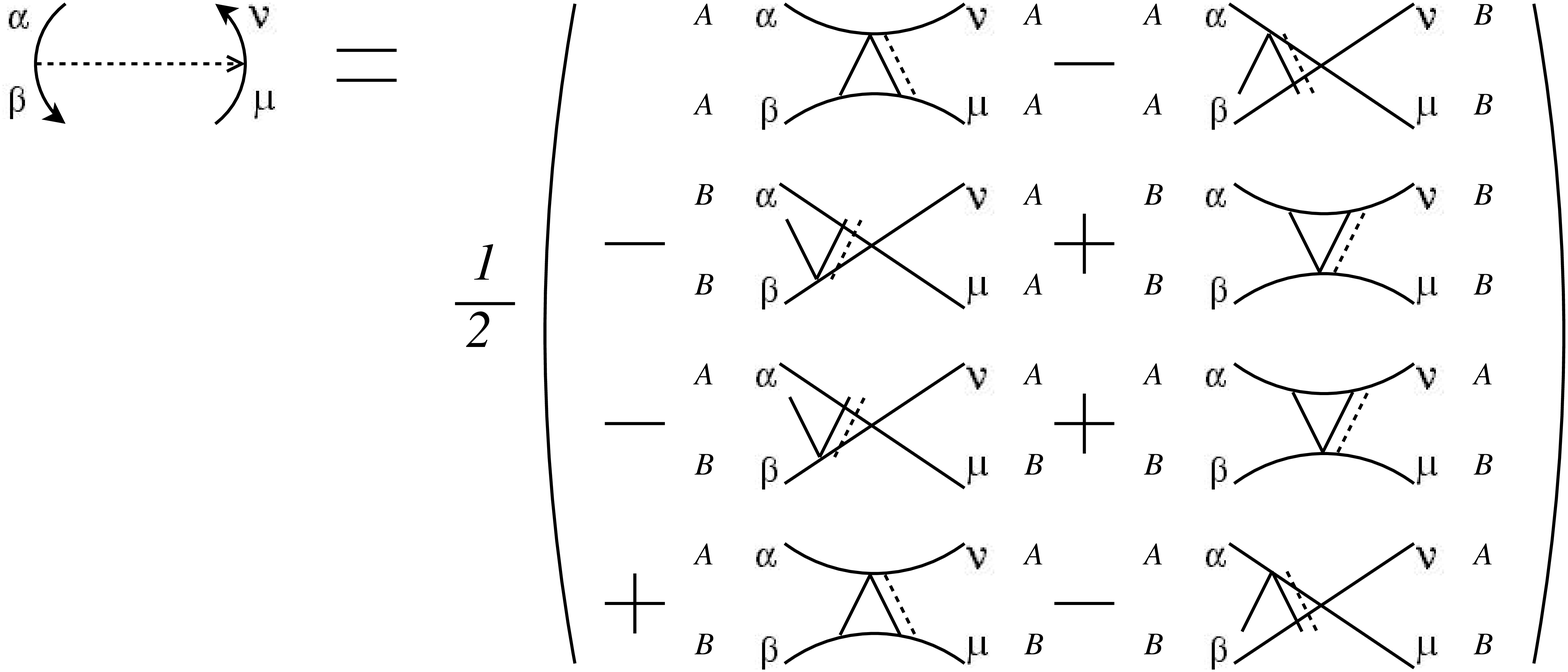, scale=0.25}}
	\caption{The $\mathfrak{so}(2N)$ tensor.}
      \label{fig:so2n}
\end{figure}

In this case the equality may be weighted with any factor because when $\alpha=\beta=\mu=\nu$ we get zero.  If we choose the factor to be $\frac{1}{2}$, as we did in the beginning of this subsection, diagramatically the $\mathfrak{so}(2N)$ $T_{\alpha\mu}^{\beta\nu}$ tensor may be expressed as in figure \ref{fig:so2n}.

\subsection{$\mathfrak{sp}(N)$}

Like in $\mathfrak{so}(2N)$ we use the basis $x_{ijk}$ and $\xi^{ijk}$ where $0\leq i\leq j\leq N$ and $k=0,1$, but the signs of some entries are different.  We have: 

\begin{eqnarray}
\nonumber x_{ij1}&=&e_{ij}-e_{j+N,i+N}\\
\nonumber x_{ij2}&=&e_{i,j+N}+e_{j,i+N}\\
\nonumber \xi^{ij1}&=&(\frac{1}{2})^{\delta^{ij}+1}(e^{ji}-e^{i+N,j+N})\\
\nonumber \xi^{ij2}&=&(\frac{1}{2})^{\delta^{ij}+1}(e^{j+N,i}+e^{i+N,j})
\end{eqnarray}

If we look at the tensor $T_{\alpha\mu}^{\beta\nu}$ we get:

\begin{eqnarray}
\nonumber T_{\alpha\mu}^{\beta\nu}&=&\sum_{i\leq j}(\xi^{ijk})_{\alpha}^{\beta}(x_{ijk})_{\mu}^{\nu}\\
\nonumber &=&\sum_{i\leq j}((\frac{1}{2})^{\delta^{ij}+1}(\delta^{\beta,j}\delta^{i}_{\alpha}-\delta^{\beta,i+N}\delta^{j+N}_{\alpha})(\delta^{\nu}_{i}\delta_{j\mu}-\delta^{\nu}_{j+N}\delta_{i+N,\mu})\\ 
\nonumber &&+(\frac{1}{2})^{\delta^{ij}+1}(\delta^{\beta,i+N}\delta^{j}_{\alpha}+\delta^{\beta,j+N}\delta^{i}_{\alpha})(\delta^{\nu}_{i}\delta_{j+N,\mu}+\delta^{\nu}_{j}\delta_{i+N,\mu}))
\end{eqnarray}

\begin{figure}[h!]
	\centerline{\epsfig{file=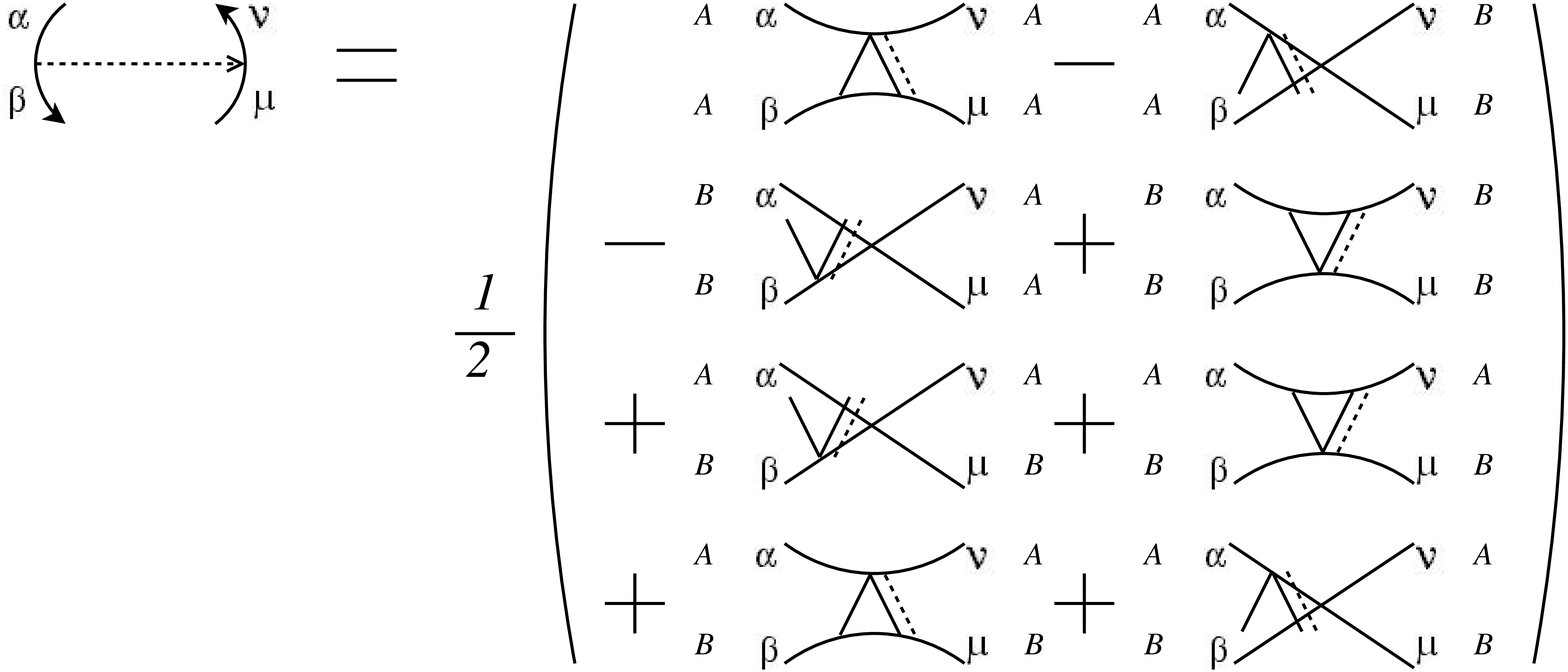, scale=0.25}}
	\caption{The $\mathfrak{sp}(N)$ tensor.}
      \label{fig:spn}
\end{figure}

The factor in the second product is due to the fact that $(x_{ii1},\xi^{ii1})=4$.  Like in the case of $\mathfrak{so}(2N)$ we assign letters $A$ and $B$ to the four corners to indicate which set each of $\alpha$, $\beta$, $\mu$ and $\nu$ belongs to.  The only possible assignments are the same as the ones for $\mathfrak{so}(2N)$.  (See figures \ref{fig:so2nABk1} and \ref{fig:so2nABk2} for the cases $k=1$ and $k=2$, respectively.)  Following the steps in the $\mathfrak{so}(2N)$ case we can see that the $\mathfrak{sp}(N)$ tensor may be given as in figure \ref{fig:spn}.  Note that it is almost the tensor for $\mathfrak{so}(2N)$ except for some sign difference.

\subsection{$\mathfrak{so}(2N+1)$}

We use the basis $x_{ijk}$ and $\xi^{ijk}$ where $0\leq i\leq j\leq N$ and $k=0,1,2$.  We have: 

\begin{eqnarray}
\nonumber x_{ij0}&=&e_{1,i+1}-e_{i+N+1,1}\\
\nonumber x_{ij1}&=&e_{i+1,j+1}-e_{j+N+1,i+N+1}\\
\nonumber x_{ij2}&=&e_{i+1,j+N+1}-e_{j+1,i+N+1}\\
\nonumber \xi^{ij0}&=&\frac{1}{2}(e_{i+1,1}-e_{1,i+N+1})\\
\nonumber \xi^{ij1}&=&(\frac{1}{2})^{\delta^{ij}+1}(e^{j+1,i+1}-e^{i+N+1,j+N+1})\\
\nonumber \xi^{ij2}&=&(\frac{1}{2})^{\delta^{ij}+1}(e^{j+N+1,i+1}-e^{i+N+1,j+1})
\end{eqnarray}

\noindent The $T_{\alpha\mu}^{\beta\nu}$ tensor is then given as

\begin{eqnarray}
\nonumber T_{\alpha\mu}^{\beta\nu}&=&\sum_{i\leq j}(\xi^{ijk})_{\alpha}^{\beta}(x_{ijk})_{\mu}^{\nu}\\
\nonumber &=&\sum_{i\leq j}(\frac{1}{2}(-\delta^{\beta,1}\delta^{i+N+1}_{\alpha}+\delta^{\beta,i+1}\delta^{1}_{\alpha})(\delta^{\nu,1}\delta^{i+1}_{\mu}-\delta^{\nu,i+N+1}\delta^{1}_{\mu})\\ 
\nonumber
&&+(\frac{1}{2})^{\delta^{ij}+1}(\delta^{\beta,j+1}\delta^{i+1}_{\alpha}-\delta^{\beta,i+N+1}\delta^{j+N+1}_{\alpha})(\delta^{\nu}_{i+1}\delta_{j+1,\mu}-\delta^{\nu}_{j+N+1}\delta_{i+N+1,\mu})\\
\nonumber &&+(\frac{1}{2})^{\delta^{ij}+1}(-\delta^{\beta,i+N+1}\delta^{j+1}_{\alpha}+\delta^{\beta,j+N+1}\delta^{i+1}_{\alpha})(\delta^{\nu}_{i+1}\delta_{j+N+1,\mu}-\delta^{\nu}_{j+1}\delta_{i+N+1,\mu}))
\end{eqnarray}

\begin{figure}[h!]
	\centerline{\epsfig{file=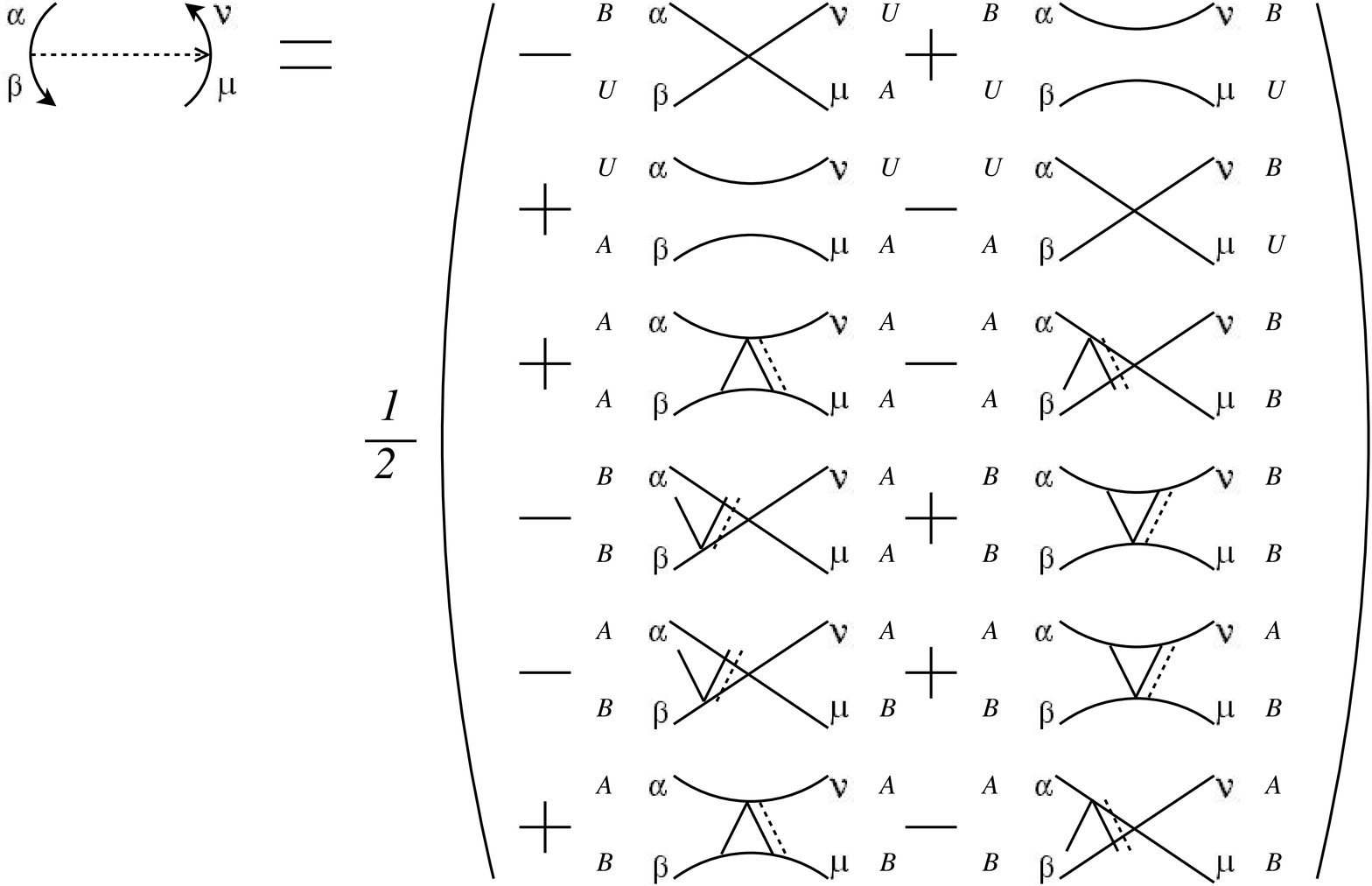, scale=0.25}}
	\caption{The $\mathfrak{so}(2N+1)$ tensor.}
      \label{fig:so2n1}
\end{figure}

Using $U$, $A$ and $B$ to denote the sets $\{1\}$, $\{2,...,N+1\}$ and $\{N+2,...,2N+1\}$ respectively, we can follow the same steps as in $\mathfrak{so}(2N)$ to get a diagrammatic representation of $T_{\alpha\mu}^{\beta\nu}$.  The result is given in figure \ref{fig:so2n1}.

\subsection{Composing the weight systems with the averaging map}

\begin{figure}[h!]
	\centerline{\epsfig{file=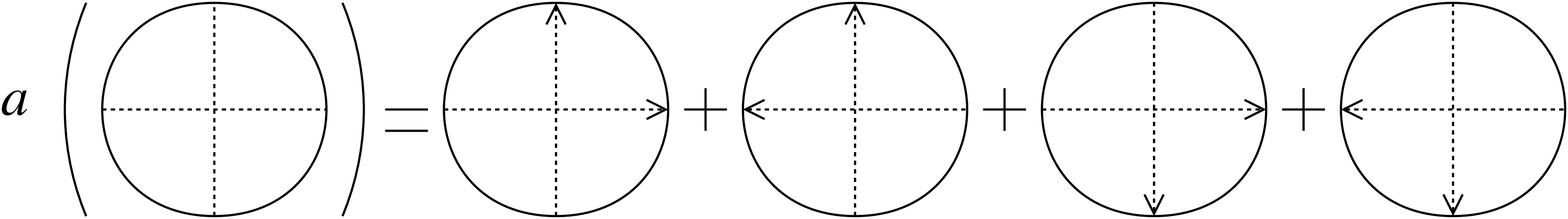, scale=0.18}}
	\caption{The averaging map.}
      \label{fig:averaging}
\end{figure}

\begin{defn}
The averaging map $a$ takes a chord diagram to an arrow diagram by summing over all possible ways to direct each chord.  (See figure \ref{fig:averaging} for an example.)
\end{defn}

Let $\mathcal{C}(\Gamma)$ be the vector space of unoriented chord diagrams on $\Gamma$ modulo the 4T relation.  (See \cite{bn1} and \cite{bn2}.)  If we have a weight system $\vec{w}:\vec{\mathcal{C}}(S^{1})\rightarrow\mathbb{C}$, then composing it with the averaging map $a$ gives us a map $w:\mathcal{C}(S^{1})\rightarrow\mathbb{C}$ which satisfies the 4T relation and is thus a weight system on $\mathcal{C}(S^{1})$.  (4T is a consequence of 6T by repeatingly applying the averaging map.)

\begin{prop}
Let $\mathfrak{g}$ be a Lie algebra in the family $\mathfrak{gl}$, $\mathfrak{sp}$ or $\mathfrak{so}$ and $(\tilde{\mathfrak{g}},\mathfrak{a},\mathfrak{a}^{*})$ be the Manin triple obtained from $\mathfrak{g}$ by the standard construction given in section \ref{sec:bialg}.  Also let $w_{\mathfrak{g}}:\mathcal{C}(S^{1})\rightarrow\mathbb{C}$ and $\vec{w}_{\tilde{\mathfrak{g}}}:\vec{\mathcal{C}}(S^{1})\rightarrow\mathbb{C}$ be the weight systems they give rise to.  Then $w_{\mathfrak{g}}=\vec{w}_{\tilde{\mathfrak{g}}}\circ a.$
\end{prop}

\begin{figure}[h!]
	\centerline{\epsfig{file=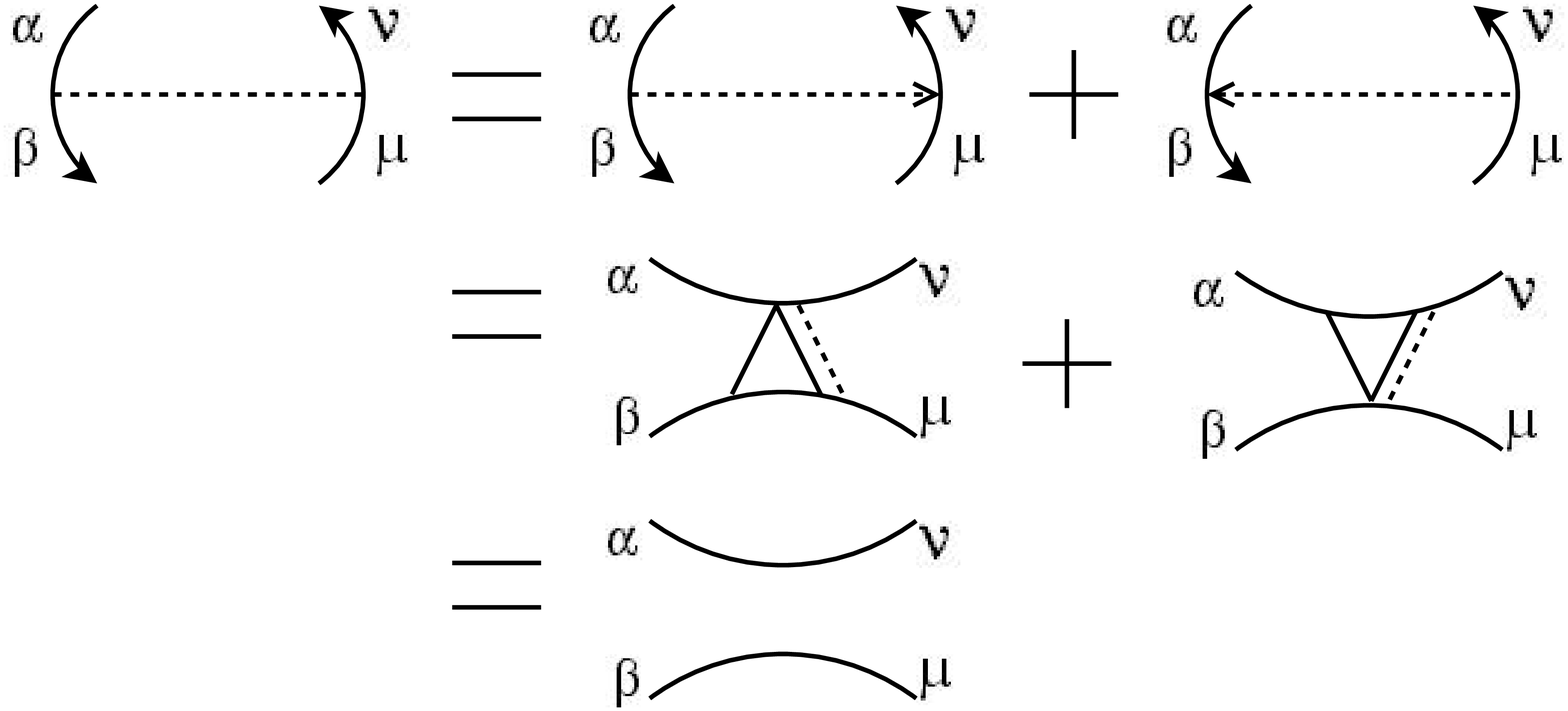, scale=0.25}}
	\caption{The unoriented $\mathfrak{gl}(N)$ tensor.}
      \label{fig:glnunori}
\end{figure}

The equality occurs at the tensor level.  We will look at each case separately.  For $\mathfrak{gl}$ the unoriented tensor is calculated as in figure \ref{fig:glnunori}.  Note the absence of restrictions on the values of unconnected Greek letters give us the weight system as given in \cite{bn1}.

\begin{figure}[h!]
	\centerline{\epsfig{file=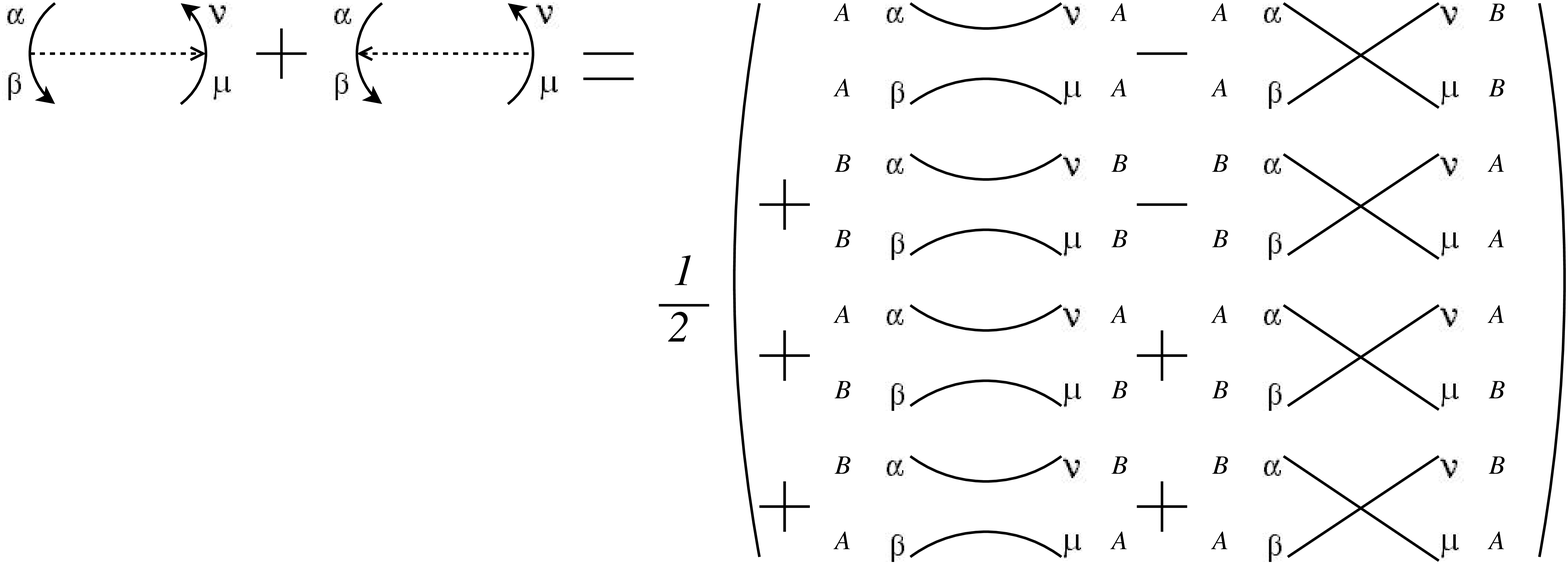, scale=0.25}}
	\caption{The unoriented $\mathfrak{sp}(N)$ tensor.}
      \label{fig:spnunori}
\end{figure}

For $\mathfrak{sp}(N)$ we have figure \ref{fig:spnunori}, which again is the weight system given in \cite{bn1}.

\begin{figure}[h!]
	\centerline{\epsfig{file=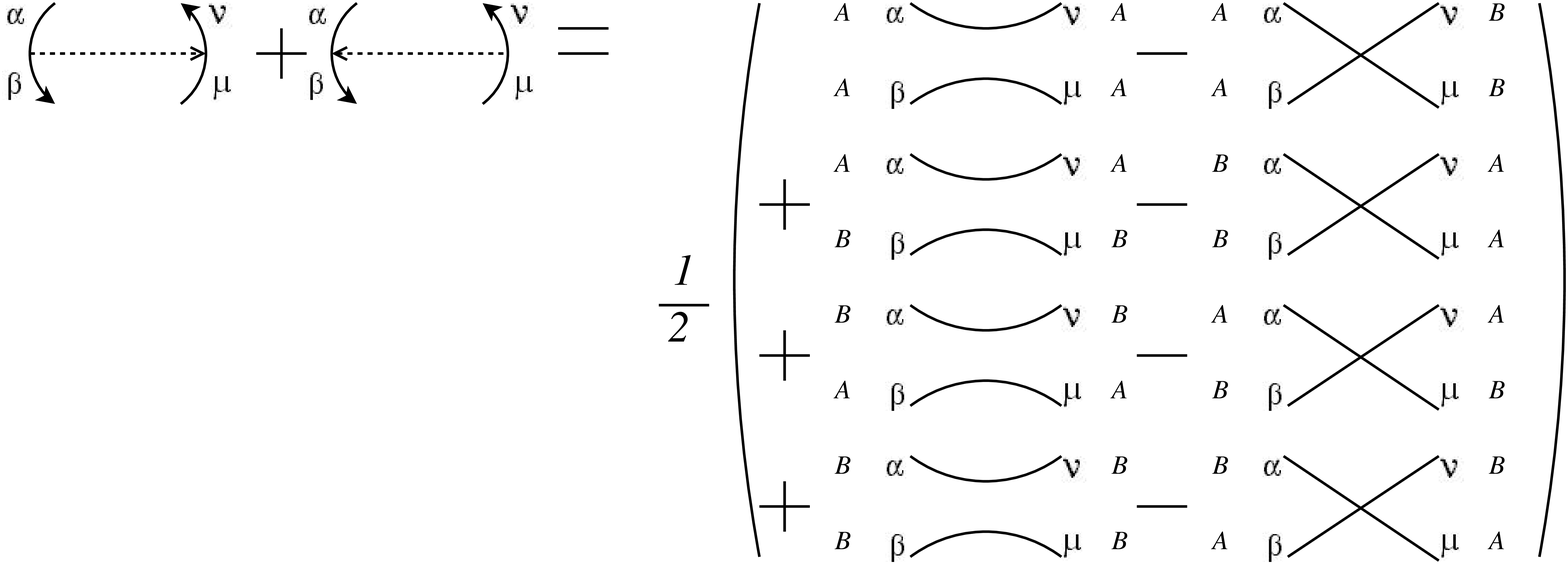, scale=0.25}}
	\caption{The unoriented $\mathfrak{so}(2N)$ tensor.}
      \label{fig:so2nunori}
\end{figure}

For $\mathfrak{so}(2N)$, we have figure \ref{fig:so2nunori}.  We consider a new basis $e'_{i}$ so that $e_{i}=M(e'_{i})$ of $\mathbb{C}^{2N}$, where 

\begin{equation}
M=\frac{1}{\sqrt{2}}\left[
\begin{array}{cc}
iI & -iI\\
-I & -I
\end{array}
\right]
\end{equation}

\noindent and $I$ is the $N\times N$ identity matrix.

\begin{figure}[h!]
	\centerline{\epsfig{file=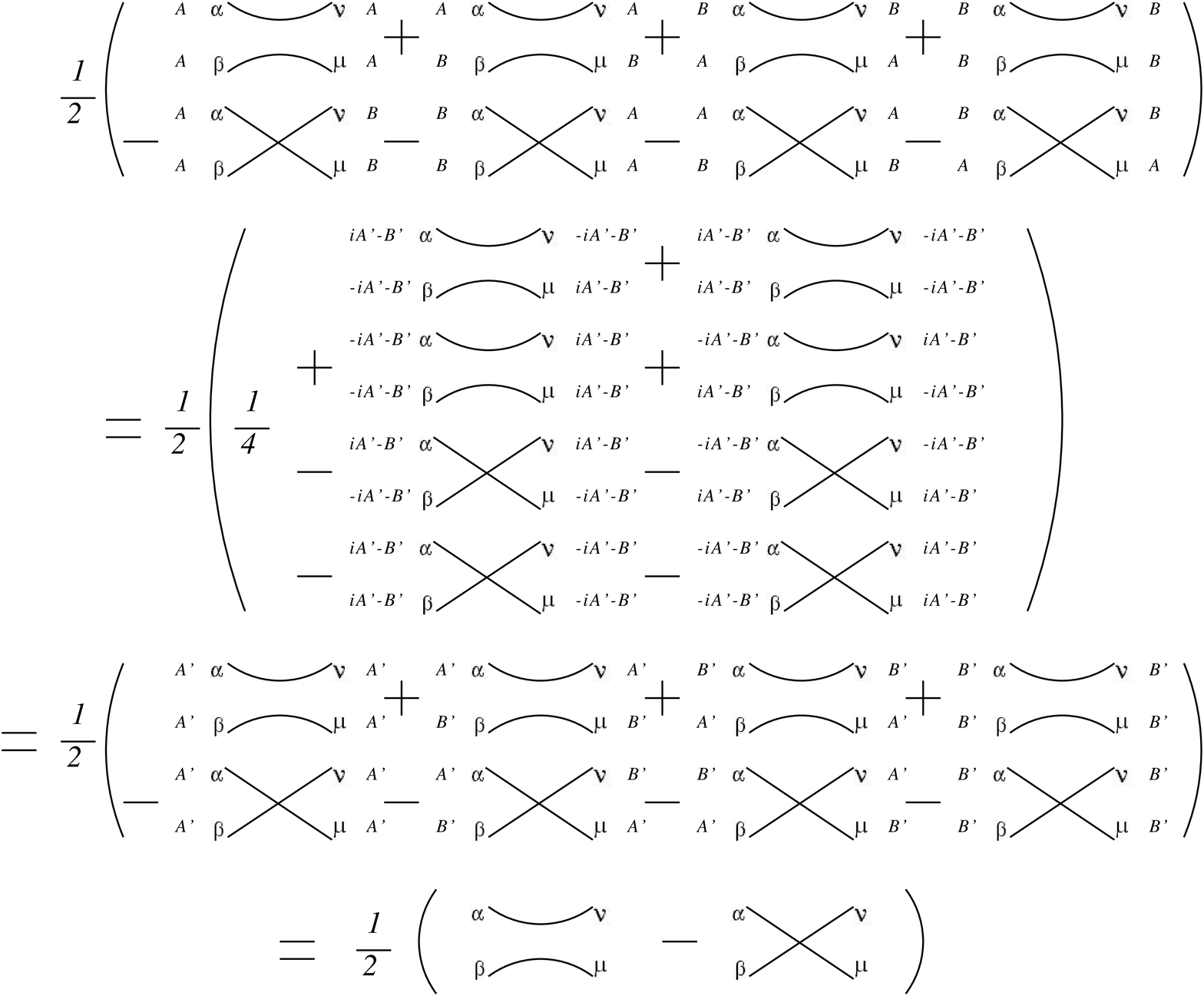, scale=0.25}}
	\caption{Calculating the unoriented $\mathfrak{so}(2N)$ tensor in the new basis.  The third line is obtained from the seond line by expansion and cancelling.  The fourth line is obtained from the third line by realizing that the latter is the former broken into cases.}
      \label{fig:so2nunorichangebasis}
\end{figure}

Note $M$ is a unitary matrix, so $(e'_{i})^{*}=\langle e'_{i},.\rangle$ where $\langle .,.\rangle$ is the inner product with respect to the basis $\{e_{i}\}$.  Using this new basis, whenever an $A$ appears at the tail of an arrow which is part of the skeleton we can replace it by $iA'-B'$.  Similarly we can replace $B$ by $-iA'-B'$.  If they appear at the head of an arrow that is part of the skeleton, however, we replace $A$ and $B$ by $-iA'-B'$ and $iA'-B'$, respectively, since the inner product is conjugate linear in the first argument.  If we expand each of the diagrams linearly we obtain the same result as in \cite{bn1}.  See figure \ref{fig:so2nunorichangebasis}.

\begin{figure}[h!]
	\centerline{\epsfig{file=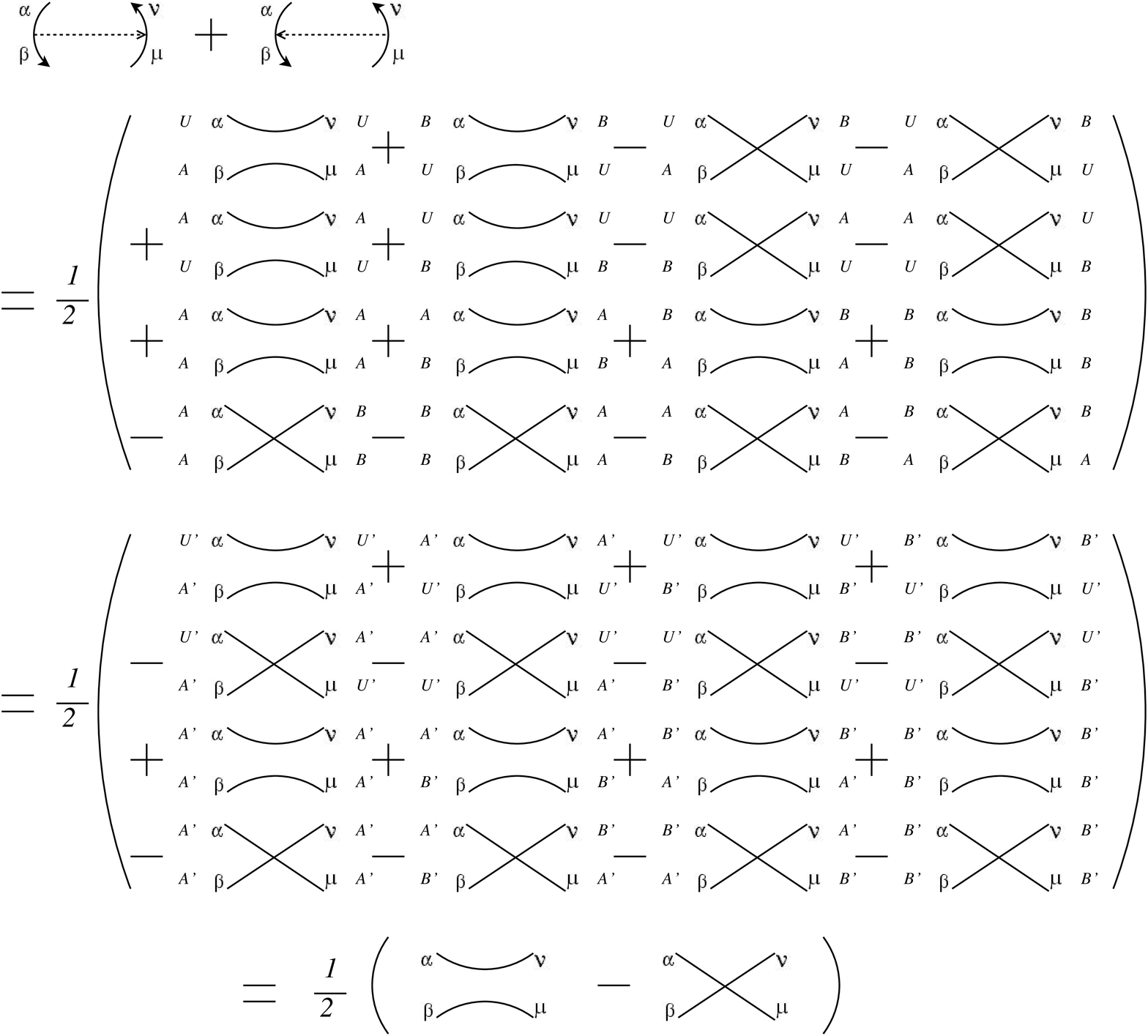, scale=0.25}}
	\caption{Calculating the unoriented $\mathfrak{so}(2N+1)$ tensor and expressing it in the new basis.}
      \label{fig:so2n1unorichangebasis}
\end{figure}

For the case $\mathfrak{so}(2N+1)$ we use the change of basis matrix

\begin{equation}
M=\frac{1}{\sqrt{2}}\left[
\begin{array}{ccc}
\sqrt{2} & 0 & 0\\
0 & iI & -iI\\
0 & -I & -I
\end{array}
\right]
\end{equation}

\noindent and we obtain the same result after expansion and cancelling.  See figure \ref{fig:so2n1unorichangebasis}.

\section{Sample calculations}
\label{sec:sample}

\begin{figure}[h!]
	\centerline{\epsfig{file=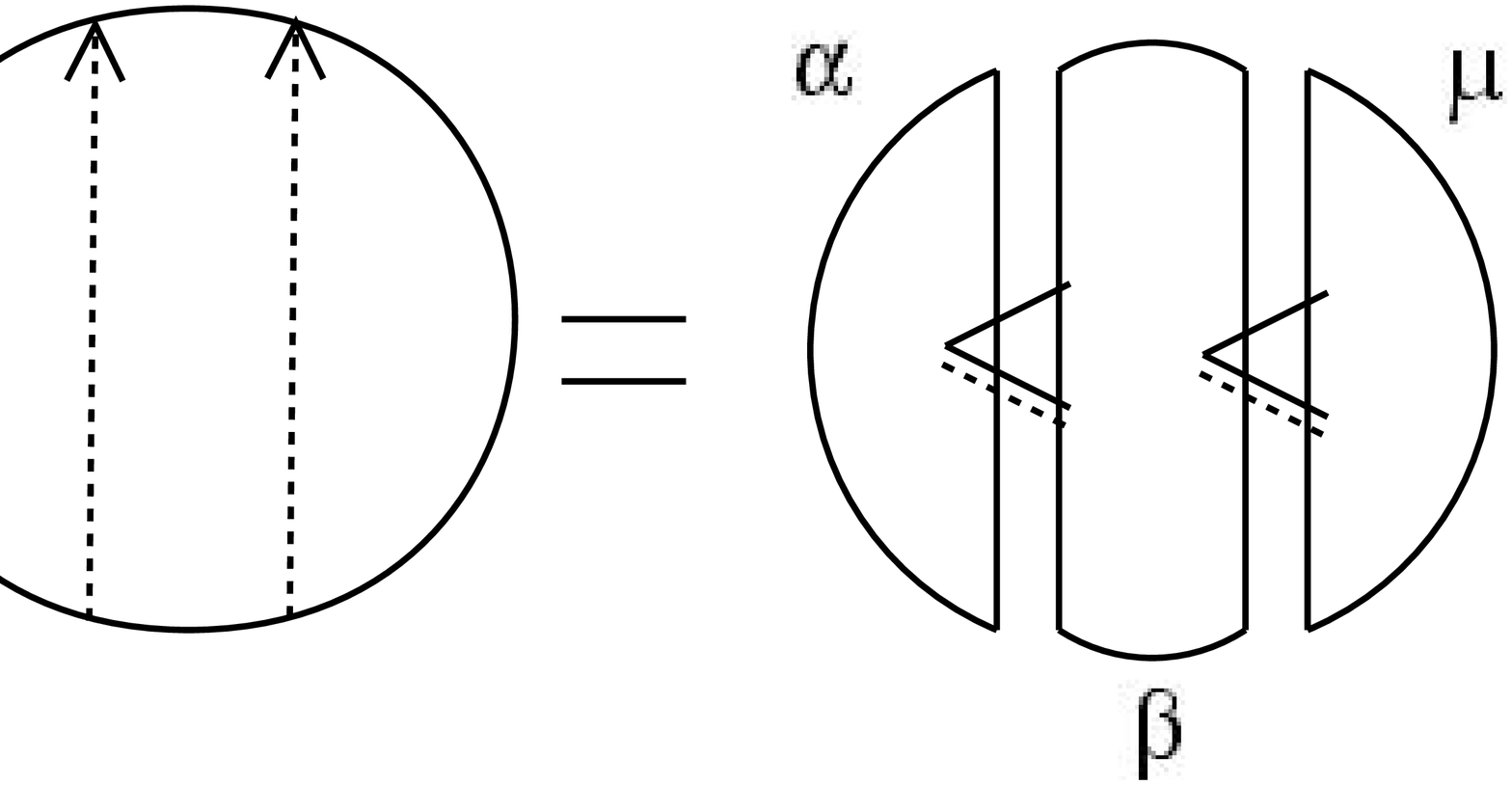, scale=0.25}}
	\caption{Calculating the $\mathfrak{gl}(N)$ weight system.}
      \label{fig:glncalc1}
\end{figure}

\begin{figure}[h!]
	\centerline{\epsfig{file=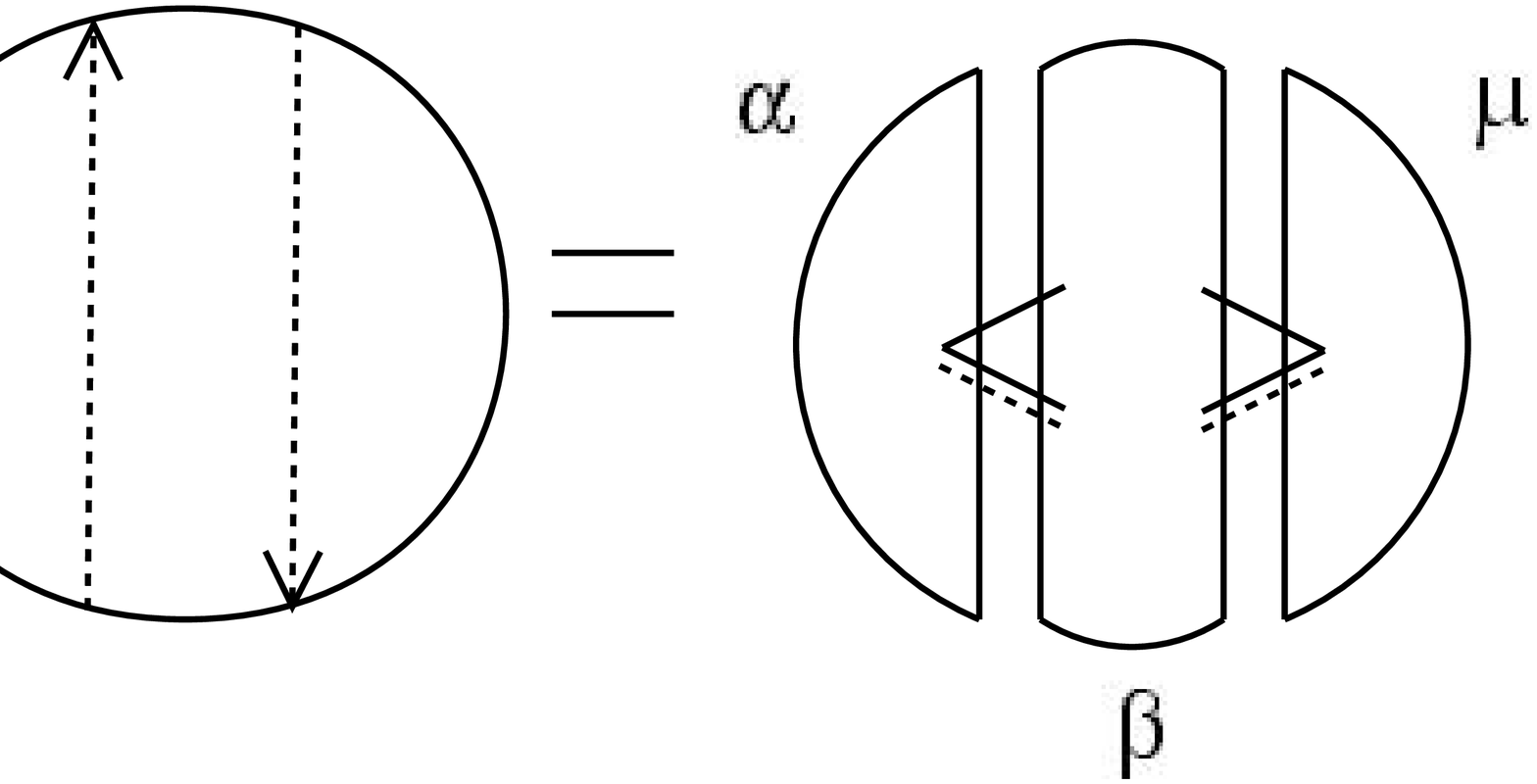, scale=0.25}}
	\caption{Calculating the $\mathfrak{gl}(N)$ weight system.}
      \label{fig:glncalc2}
\end{figure}

We now do some sample calculations.  In this section the skeleton is always oriented counterclockwise.  First we calculate the two diagrams shown in figures \ref{fig:glncalc1} and \ref{fig:glncalc2}, using the $\mathfrak{gl}(N)$ weight system.  For the first picture, each triple $(\alpha, \beta, \mu)\in\{1,...,N\}^{3}$ such that $\alpha\leq\beta\leq\mu$ gives us a value of either 1, $\frac{1}{2}$ or $\frac{1}{4}$, depending on whether $\alpha=\beta$ or $\beta=\mu$.  Therefore the diagram should have weight $a+\frac{1}{2}b+\frac{1}{4}c$, where $a$ is the number of triples $(\alpha, \beta, \mu)$ such that $\alpha<\beta<\mu$, $b$ is the number of triples $(\alpha, \beta, \mu)$ such that $\alpha<\beta=\mu$ or $\alpha=\beta<\mu$, and $c$ is the number of triples $(\alpha, \beta, \mu)$ such that $\alpha=\beta=\mu$.  The number is

\begin{center}
${N\choose 3}+{N\choose 2}+\frac{N}{4}$.
\end{center}

\noindent Using a similar argument we know the weight of the picture in figure \ref{fig:glncalc2} is

\begin{center}
$2{N\choose 3}+2{N\choose 2}+\frac{N}{4}$,
\end{center}

\noindent so the $\mathfrak{gl}(N)$ weight system is capable of telling the two diagrams apart.

\begin{figure}[h!]
	\centerline{\epsfig{file=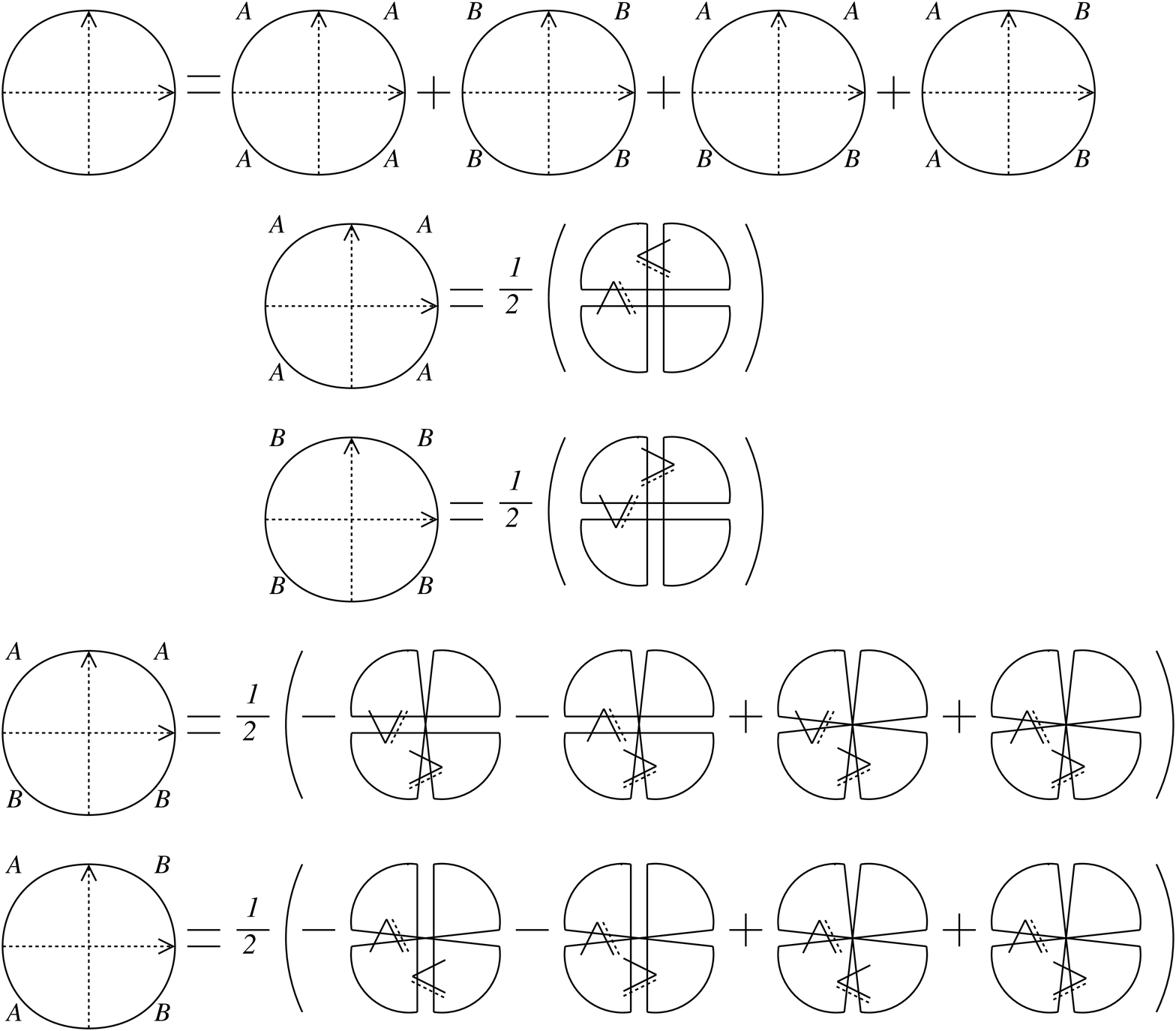, scale=0.20}}
	\caption{Calculating the $\mathfrak{so}(2N)$ weight system.}
      \label{fig:so2ncalc}
\end{figure}

Now we calculate the weight of the picture in figure \ref{fig:so2ncalc} using $\mathfrak{so}(2N)$.  We assign letters $A,B$ to each arc following the rules from section \ref{sec:calc} (see figure \ref{fig:so2ncalc}).  For each assignment we have one or four ways to resolve the diagram, and for a resolution with $k$ loops we count the number of $k$-tuples in $\{1,...,N\}$ such that the inequalities are satisfied, bearing in mind that each equality comes with a weight $\frac{1}{2}$.  The weights of the first two diagrams and the last two diagrams on the right hand side of the first equation in figure \ref{fig:so2ncalc} are therefore $\frac{1}{2}(\frac{N}{4})$ and $\frac{1}{2}(-\frac{N}{4}-\frac{N}{4}+\frac{N}{4}+({N\choose 2}+\frac{N}{4}))$, respectively.  The weight of the diagram is therefore

\begin{center}
$\frac{N}{4}+{N\choose 2}$.
\end{center} 

\begin{figure}[h!]
	\centerline{\epsfig{file=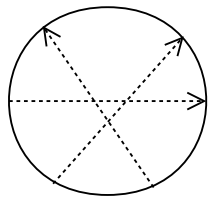, scale=0.25}}
	\caption{A diagram with three arrows.}
      \label{fig:threearrows}
\end{figure}

\begin{figure}[h!]
	\centerline{\epsfig{file=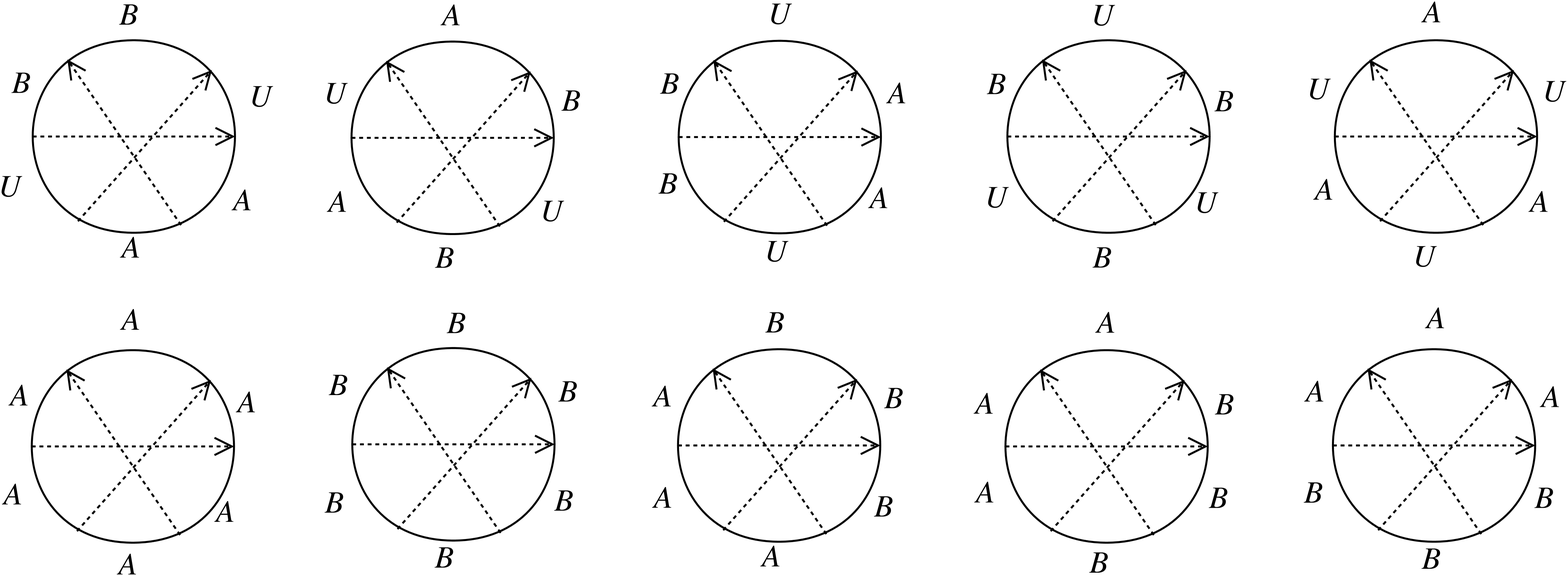, scale=0.20}}
	\caption{The ten ways to label the arcs of the skeleton for the calculation of the $\mathfrak{so}(2N+1)$ weight system.}
      \label{fig:so2n1samplelabel}
\end{figure}

Finally we calculate the weight of the picture in figure \ref{fig:threearrows}.  According to section \ref{sec:calc} we have ten possible ways to label the arcs on the skeleton (figure \ref{fig:so2n1samplelabel}).  The weights of the ten diagrams (first row left to right, then second row left to right) are $-{N\choose 2}-\frac{N}{2}$, $-2{N\choose 2}$, $-{N\choose 2}-\frac{N}{2}$, $N$, $N$, $\frac{N}{8}$, $\frac{N}{8}$, $-{N\choose 3}-{N\choose 2}$, $-2{N\choose 3}-{N\choose 2}-\frac{N}{8}$, and $-{N\choose 3}-{N\choose 2}+\frac{N}{4}$, respectively.  Summing these ten weights we get that the weight of figure \ref{fig:threearrows} is $-4{N\choose 3}-7{N\choose 2}+\frac{11N}{8}$.

\pagebreak


\begin{thebibliography}{99}

\bibitem[BN1]{bn1}D. Bar-Natan.  \textit{Weights of Feynman Diagrams and the Vassiliev Knot Invariants}.  http://www.math.toronto.edu/~drorbn/papers/weights/weights.ps

\bibitem[BN2]{bn2}D. Bar-Natan.  \textit{On the Vassiliev Knot Invariants}.  Topology \textbf{34} (1995), 423-472.

\bibitem[CP]{cp}V. Chari and A. Pressley.  \textit{A Guide to Quantum Groups}.  Cambridge University Press, 1994

\bibitem[ES]{es}E. Etingof and O. Schiffman. \textit{Lectures on Quantum Groups}. International Press, Boston, 1998.  

\bibitem[GPV]{gpv}M. Goussarov, M. Polyak and O. Viro. \textit{Finite Type Invariants of Classical and Virtual Knots}. Topology \textbf{39} (2000), no. 5, 1045-1068.  arXiv:math.GT/9810073


\bibitem[Ha]{haviv}A. Haviv. \textit{Towards a Diagrammatic Analogue of the Reshetikhin-Turaev Link Invariants}. 2002.  arXiv.math.QA/0211031

\bibitem[Po]{polyak}M. Polyak. \textit{On the Algebra of Arrow Diagrams}. Letters in Mathematical Physics \textbf{51} (2000): 275-291. 




\end{thebibliography}
\end{document}